\documentclass[11pt]{article}

\usepackage{amsfonts,amsmath,amsthm,amssymb}
\usepackage{bm}
\usepackage[colorlinks, citecolor=blue, urlcolor=blue, backref=page]{hyperref}

\usepackage{geometry}
\usepackage[utf8]{inputenc}
\usepackage[T1]{fontenc}    % use 8-bit T1 fonts

\usepackage{enumitem}
\usepackage{hyperref}
\usepackage{url}

\newtheorem{theorem}{Theorem}[section]

\newtheorem{assumption}[theorem]{Assumption}

\theoremstyle{definition}

\newtheorem{example}[theorem]{Example}
\newtheorem{remark}[theorem]{\textbf{Remark}}
\numberwithin{equation}{section}

\newcommand{\e}{{\rm e}}

\renewcommand{\P}{{\mathbb P}}

\newcommand{\E}{{\mathbb E}}

\newcommand{\Q}{{\mathbb Q}}
\newcommand{\R}{{\mathbb R}}

\newcommand{\N}{{\mathbb N}}

\newcommand{\Fcal}{{\mathcal F}}

\newcommand{\Ical}{{\mathcal I}}

\newcommand{\Mcal}{{\mathcal M}}

\newcommand{\Tcal}{{\mathcal T}}

\newcommand{\ir}{{i}}
\newcommand{\jr}{{j}}

\newcommand{\nG}{{G}}

\numberwithin{equation}{section}
\numberwithin{theorem}{section}

\begin{document}

\title{Propagation of chaos for point processes induced by particle systems with mean-field drift interaction\footnote{This work has been partially supported by the National Science Foundation under grant NSF DMS-2206062.}}

\author{Nikolaos Kolliopoulos\thanks{Department of Mathematics, University of Michigan, \url{nkolliop@andrew.cmu.edu}.}
\and
Martin Larsson\thanks{Department of Mathematical Sciences, Carnegie Mellon University, \url{larsson@cmu.edu}.}
\and
Zeyu Zhang\thanks{Department of Mathematical Sciences, Carnegie Mellon University, \url{feketerigoem@gmail.com}.}
}
%\date{}
\maketitle

\begin{abstract}
% We studied the asymptotic behavior of the normalized maxima of real-valued diffusive particles with mean-field drift interaction in our previous work and established the result of propagation of chaos. This is a subsequential work study about the asymptotic behavior of the point process induced by the same interacting particle system, which contains the information of the joint distribution of arbitrary top $k$ order statistics. We manager to establish the propagation of chaos for such point process: it has the same weak limit as the point process induced by the limit Mckean-Vlasov dynamics. The later is exactly Poisson Point Process if certain assumptions holds.

%We studied the asymptotic behavior of the normalized maxima of real-valued diffusive particles with mean-field drift interaction in our previous work and established the result of propagation of chaos. This is a subsequential work study about the asymptotic behavior of the point process induced by the same interacting particle system. Under suitable assumptions, we manage to establish propagation of chaos for this point process: it has the same weak limit as the point process induced by i.i.d. copies of the solution of a limiting Mckean--Vlasov equation. This weak limit is a Poisson point process whose mean measure is related to classical extreme value distributions.

We study the asymptotics of the point process induced by an interacting particle system with mean-field drift interaction. Under suitable assumptions, we establish propagation of chaos for this point process: it has the same weak limit as the point process induced by i.i.d. copies of the solution of a limiting Mckean--Vlasov equation. This weak limit is a Poisson point process whose intensity measure is related to classical extreme value distributions. In particular, this yields the limiting distribution of the normalized upper order statistics.
\end{abstract}

%\tableofcontents

\section{Introduction}

%concerned with the large-population asymptotics of the maxima of certain real-valued diffusive particle systems $X^{1,N},\ldots,X^{N,N}$ with mean-field interaction through the drifts. Specifically, we are interested in large-$N$ limits of}

This paper is concerned with the upper edge asymptotics of a system of real-valued diffusive particles $X^{1,N},\ldots,X^{N,N}$ with mean-field interaction through the drifts, as the population size $N$ tends to infinity. This is done through a study of
the point process $\widetilde U_N$ on $[0,1]\times\mathbb{R}$ defined by
\begin{equation}\label{pointprocess}
\widetilde{U}_N(I) = \sum_{i=1}^N \bm 1_{\left\{\left(\frac{i}{N}, (X^{i,N}_T - b_T^N) / a_T^N \right)\in I\right\}}
%U_N(I) = \sum_{i=1}^N \bm 1_{\left\{\left(\frac{i}{N}, \frac{X^{i,N}_T - b_T^N}{a_T^N}\right)\in I\right\}},
\end{equation}
for any Borel subset $I \subset [0,1]\times\R$, where $a^N_T$ and $b^N_T$ are suitable normalizing constants. In particular, this allows us to derive the large $N$ asymptotics of any fixed number $k$ of upper order statistics. Our result can be viewed as a propagation of chaos property for the point process $\widetilde U_N$. Alternatively, it can be viewed as an extension of classical extreme value theory for i.i.d.\ random variables to the case of mean-field interaction.
This addresses a natural question left open in our earlier work \cite{kolliopoulos2022propagation}, which establishes propagation of chaos for the largest particle only.

% This paper extends our earlier work \cite{kolliopoulos2022propagation}, which establishes propagation of chaos for the largest particle in a system of real-valued diffusive particles $X^{1,N},\ldots,X^{N,N}$ with mean-field interaction through the drifts. In the present paper we do not only consider the maximum particle, but we study the point process $\widetilde U_N$ on $[0,1]\times\mathbb{R}$ defined by
% \begin{equation}\label{pointprocess}
% \widetilde{U}_N(I) = \sum_{i=1}^N \bm 1_{\left\{\left(\frac{i}{N}, (X^{i,N}_T - b_T^N) / a_T^N \right)\in I\right\}}
% %U_N(I) = \sum_{i=1}^N \bm 1_{\left\{\left(\frac{i}{N}, \frac{X^{i,N}_T - b_T^N}{a_T^N}\right)\in I\right\}},
% \end{equation}
% for any Borel subset $I \subset [0,1]\times\R$, where $a^N_T$ and $b^N_T$ are suitable normalizing constants. In particular, this allows us to derive the large $N$ asymptotics for any fixed number $k$ of upper order statistics. %, since the point process contains all the information about the finite joint marginal distribution of the normalized ranked system at fixed time point. 
% Our result can be viewed as a propagation of chaos property for the point process $U_N$. Alternatively, it can be viewed as an extension of classical extreme value theory for i.i.d.\ random variables to the case of mean-field interaction.

%This is an extension of the propagation of chaos property for the point process induced by particle systems, as well as an extension of the classical extreme value theory to interacting situation than the i.i.d one.

We work with the same particle dynamics as in \cite{kolliopoulos2022propagation}. Specifically, for each $N \in \N$ the $N$-particle system evolves according to a stochastic differential equation of the form
\begin{equation} \label{eq_particle_SDE}
dX^{i,N}_t = A(t, X^{i,N}_{[0,t]}) \left( B\left( t, X^{i,N}_{[0,t]}, \int g(t, X^{i,N}_{[0,t]}, \bm y_{[0,t]}) \mu^N_t(d\bm y) \right) dt + dW^i_t \right) + C(t, X^{i,N}_{[0,t]}) dt
\end{equation}
for $i=1,\ldots,N$, with i.i.d.\ initial conditions $X^{i,N}_0 \sim \nu_0$, where $\nu_0$ is a given probability measure on $\R$. Moreover, $W^i$, $i \in \N$, is family of independent standard Brownian motions, and we write $\bm x_{[0,t]} = (x(s))_{s \in [0,t]}$ for any continuous function $\bm x$. For each $t \in \R_+$ we defined the following empirical measure on trajectory up to time t,
\begin{equation}\label{iem}
\mu^N_t = \frac1N \sum_{i=1}^N \delta_{X^{i,N}_{[0,t]}}.
\end{equation}
The coefficients $A(t, \bm x_{[0,t]})$, $B(t, \bm x_{[0,t]}, r)$, $C(t, \bm x_{[0,t]})$ and interaction function $g(t, \bm x_{[0,t]}, \bm y_{[0,t]})$ are defined for all $t \in \R_+$, $\bm x, \bm y \in C(\R_+)$, and $r \in \R$. Note that here we only allow interaction in the drift $B$ but not in $A$ or $C$. This restriction is crucial for the methods used later, which rely on the ability to use equivalent changes of probability measure to remove the interaction.

Systems that fall directly into the setup of this paper have been used in physics \cite{JAWA18, SY20}, in neuroscience \cite[Display (1.2)]{LP19}, and also in the modelling of the parameters of learning methods (\cite[Display (52)]{MMM19} and \cite[Display (114)]{ROE18}). The credit risk mean-field game model studied in \cite{CFS15} is also of the form we consider here under Nash equilibrium, up to the direct dependence of the equilibria on $1/N$ which vanishes as $N \to \infty$ and is not expected to lead to any major difficulty. The systems used in \cite{BOCA16, BOPARO18} for the modelling of interbanking networks differ from a special case of our setup (see example~\ref{ex_1}) only in the presence of a common jump term, so the method we develop here could reduce the study of their upper order statistics to that of simpler systems with no mean-field interaction. In Stochastic Portfolio Theory, a study of the upper order statistics is required to check whether the empirically observed linearity of the upper-left edge of the capital distribution curve is captured by certain mean-field models like the Atlas model \cite{REY17}; even though mean-field interaction in the Atlas model is of the form studied in this paper, it is also present in the noise terms, and further extension of our work is required (see example~\ref{ex_2} for a variation with interaction only in the drifts).

\color{black} Classical propagation of chaos \cite{MR221595,MR1108185,MR968996, MR3841406} states, under suitable assumption, that for any fixed number $k \in \N$, the first $k$ particles $(X^{1,N}, \ldots, X^{k,N})$ converge jointly as $N \to \infty$ to $k$ independent copies $(X^1,\ldots,X^k)$ of the solution to the McKean--Vlasov equation
\begin{equation} \label{eq_McKV_SDE}
\begin{aligned}
    dX_t &= A(t, X_{[0,t]}) \left( B\left( t, X_{[0,t]}, \int g(t, X_{[0,t]}, \bm y_{[0,t]}) \mu_t(d\bm y) \right) dt + dW_t \right) + C(t, X_{[0,t]}) dt \\
\mu_t &= \text{Law}(X_{[0,t]})
\end{aligned}
\end{equation}
with initial condition $\mu_0 = \nu_0$. See in particular \cite[Theorem~2.1]{JA19}. This leads to the heuristic that the interacting $N$-particle system behaves asymptotically like a system of i.i.d.\ copies of the solution of the limiting McKean--Vlasov equation. This suggests that the large-$N$ asymptotics of the point process $\widetilde{U}_N$ in \eqref{pointprocess} should match the asymptotics of the point process $U_N$ induced by the independent copies $X^i$ of the solution of \eqref{eq_McKV_SDE}, namely
\begin{equation} \label{eq_point_process_McKV}
U_N(I) = \sum_{i=1}^N \bm1_{\left\{\left(\frac{i}{N}, (X^{i}_T - b_T^N) / a_T^N\right)\in I\right\}}.
\end{equation}
%
%From the perspective of the classical propagation of chaos: one would expect that the large population of the interacting particle system behaves asymptotically like a system of i.i.d copies of the solution of limiting McKean-Vlasov equation. From the perspective of extreme value theory: for $N$-i.i.d copies of any random variable, the order statistics of them should converge weakly to some known distribution (which we call it as domain of attraction) and the point process induced by them converge weakly to Poisson Point process.  
%
%Such heuristic suggests that the large-$N$ asymptotic of the point process in \eqref{pointprocess} should match the asymptotic of the point process induced by the independent copies $X^i$ of the solution of \eqref{eq_McKV_SDE}:
%\begin{equation} \label{eq_point_process_McKV}
%\tilde{U}_N(I) = \sum_{i=1}^N \bm1_{\left\{\left(\frac{i}{N}, \frac{X^{i}_T - b_T^N}{a_T^N}\right)\in I\right\}}.
%\end{equation}
Due to the i.i.d.\ property, the latter falls within the framework of classical extreme value theory, and its asymptotic behavior is therefore well-known. Specifically, under suitable assumptions, $U_N$ converges weakly to a Poisson random measure whose intensity measure depends on the so-called extreme value index of the distribution of $X^1_T$. See \cite[Theorem 2.1.2]{HF10} for a precise statement and, more generally, \cite{HF10, MR2364939} for an introduction about extreme value theory.

Our main result, Theorem~\ref{T_main_precise}, states (under assumptions) that the above intuition is correct: $\widetilde{U}_N$ converges to the same Poisson point process as $U_N$. However, to prove this, existing results on standard propagation of chaos, such as the rate of convergence statements in \cite{JA19,LA22}, appear to be insufficient. Instead, we rely on an extension of a method developed in our earlier paper \cite{kolliopoulos2022propagation}. As in \cite{dawsont1987large,MR3841406,JA19}, our method starts with a Girsanov transformation. We then depart from earlier methods by performing an iterated integral expansion of the Radon--Nikodym derivative. By carefully controlling the terms in this expansion, we are able to complete the proof.

%However, the above intuition is flawed. From the inequality in Theorem 2.1 in \cite{JA19}, one may not expect the existence of the weak limit of the following form:
%\[
%f(X_t^{1,N},X_t^{2,N},\cdots, X_t^{k_N,N}).
%\]
%
%Apart from this, \cite{LA22} result indicate that the relative entropy between law of 
%\[
%(X^{1,N}_t,\cdots,X^{k,N}_t)
%\] and 
%\[
%(X^1_t,\cdots,X^{k}_t)
%\]
%has an convergence rate of $(k/N)^2$. Both suggests the standard argument of propagation of chaos fall in our situation. We have established an method in \cite{kolliopoulos2022propagation22} used measure change followed by an detailed analysis of the iterated integrals appeal in chaos-expansion of the Radon–Nikodym density to proof the propagation of chaos for the maxima of the interacting particle system. Following the similar spirit, in this paper, we extend our results to establish the propagation of chaos for the Point process defined in \eqref{pointprocess}, which is precisely the main theorem \ref{T_main_precise} . Consequently, the result of propagation of chaos for the order statistics along with their joint distribution have been established.

The rest of the paper organized as follow: we discuss the assumptions, notation, and state our main result precisely in Section~\ref{S_results}. We then discuss some examples and consequences of our main result in Section~\ref{S_examples}. In Section~\ref{S_main_proof}, we give the proof of the main result, Theorem~\ref{T_main_precise}. We will frequently use the notation
\[
[n] = \{1,\ldots,n\} \text{ for any } n \in \N = \{1,2,\ldots\},
\]
and $\R_+ = [0,\infty)$. We will allow generic constants $C$ to vary from line to line, and occasionally indicate the dependence on parameters by writing $C(n)$, $C(p,n)$, etc. We also use the notation
\[
\Fcal^V_t = \sigma(X^i_s,W^i_s \colon s \le t, i \in V) \text{ for any index set } V \subset \N,
\]
and write $\mathbb L$ for the space of all progressively measurable processes $Y$ with locally integrable moments.

\section{Assumptions and main result} \label{S_results}

We work in the same setup as in our previous work, and under the same assumption; see \cite[Assumption 2.1, Assumption 2.3]{kolliopoulos2022propagation}. For the convenience of the reader, we recall these here.

\begin{assumption}\label{ass1}
The coefficient functions $(t,\bm x) \mapsto A(t, \bm x_{[0,t]})$, $(t,\bm x,r) \mapsto B(t, \bm x_{[0,t]}, r)$, $(t,\bm x) \mapsto C(t, \bm x_{[0,t]})$ and the interaction function $(t,\bm x,\bm y) \mapsto g(t, \bm x_{[0,t]}, \bm y_{[0,t]})$
are real-valued measurable functions on $\R_+ \times C(\R_+)$, $\R_+ \times C(\R_+) \times \R$, $\R_+ \times C(\R_+)$, and $\R_+ \times C(\R_+) \times C(\R_+)$, respectively. They satisfy the following conditions:
\begin{itemize}
\item $A$ and $C$ are uniformly bounded,
\item for every $t \in \R_+$ and $\bm x \in C(\R_+)$, the function $r \mapsto B(t, \bm x_{[0,t]}, r)$ is twice continuously differentiable, and its first and second derivatives are bounded uniformly in $(t, \bm x)$.
\end{itemize}
\end{assumption}

As in \cite{kolliopoulos2022propagation}, we avoid imposing specific conditions for well-posedness of the limiting equation, but rather assume existence directly (uniqueness is not actually required, so we do not assume it.) We also assume certain moment bounds. Two examples where these assumptions can be verified are discussed below; see Examples~\ref{ex_1} and~\ref{ex_2}. See also \cite[Remark~2.7 and Section~3]{kolliopoulos2022propagation}.

%By imposing further conditions we could appeal to known results on well-posedness of McKean--Vlasov equations to assert that \eqref{eq_McKV_SDE} has a solution. Rather than doing this, we will assume existence directly (uniqueness is not actually required, so we do not assume it.) 

\begin{assumption}\label{ass2}
Fix a probability measure $\nu_0$ on $\R$ and assume that the McKean--Vlasov equation \eqref{eq_McKV_SDE} admits a weak solution $(X,W)$ with $X_0 \sim \nu_0$. Construct (for instance as a countable product) a filtered probability space $(\Omega,\Fcal,(\Fcal_t)_{t\ge0},\P)$ with a countable sequence $(X^i,W^i)$, $i \in \N$ of independent copies of $(X,W)$. Then, assume that there is a continuous function $K(t)$ such that for all $p \in \N$, $t \in \R_+$, $N \in \N$, and $i,j \in [N]$, one has the moment bounds
\begin{equation}\label{part1ass2}
\E\left[ g(t, X^i_{[0,t]}, X^j_{[0,t]})^{2p} \right] \le p!\, K(t)^p
\end{equation}
and
\begin{equation}\label{part2ass2}
\E\left[ \left( \int g(t, X^i_{[0,t]}, \bm y_{[0,t]}) ( \mu^N_t - \mu_t)(d\bm y) \right)^{2p} \right] \le \frac{1}{N^p} p! \, K(t)^p.
\end{equation}
\end{assumption}

%One may refer to remark 2.7 in [previous work] for further discussion about the sufficient conditions for the moment bounds \eqref{part1ass2}--\eqref{part2ass2} and the well-posedness of McKean-Vlasov equations.
We used the same notation as in \cite{MR2364939}. In general a point process $U$ is a random distribution of points in a state space $E$ which is locally compact with a countable basis. Let $\varepsilon$ be the Borel $\sigma$-algebra of subsets of $E$. More precisely letting $M_p(E)$ denote the space of all point measures on $E$ and define a $\sigma$-algebra $\Mcal_p(E)$ of subsets of $M_p(E)$ to be the smallest $\sigma$-algebra making all evaluation maps $m\to m(F)$ measurable for all $F\in\varepsilon$. Then $U$ is a measurable map from some probability space $(\Omega, \mathcal{F})$ into $(M_p(E), \mathcal{M}_p(E))$. Readers may refer to \cite[Section 3]{MR2364939} for the background of point process. In our situation, the state space $E$ for $U_N, \widetilde{U}_N$ is $[0,1]\times \mathbb{R}$.

In preparation for the Theorem~\ref{T_main_precise}, we now introduce some notation and definitions. % from our previous work \cite{kolliopoulos2022propagation22}.
Suppose Assumptions~\ref{ass1} and~\ref{ass2} are satisfied. For each $N \in \N$ we use the processes $(X^i,W^i)$ to construct the $N$-particle systems through a change of probability measure. The $N$-particle empirical measure is defined by
\begin{equation} \label{eq_muN}
\mu^N_t = \frac1N \sum_{i=1}^N \delta_{X^i_{[0,t]}}.
\end{equation}
We then define a Radon--Nikodym density process by
\begin{equation}\label{RND}
Z^N = \exp(M^N - \frac12\langle M^N \rangle)
\end{equation}
where
\begin{equation} \label{eq_MN}
M^N_t = \sum_{i=1}^N \int_0^t \Delta B^{i,N}_s dW^i_s
\end{equation}
and
\begin{equation} \label{eq_DeltaBN}
\begin{aligned}
\Delta B^{i,N}_t &= B\left(t,X_{[0,t]}^i,\int g(t, X_{[0,t]}^i, \bm{y}_{[0,t]}) \mu_t^N(d\bm{y})\right) \\
& \qquad - B\left(t,X_{[0,t]}^i,\int g(t, X_{[0,t]}^i, \bm{y}_{[0,t]}) \mu_t(d\bm{y})\right).
\end{aligned}
\end{equation}
Under Assumptions~\ref{ass1} and~\ref{ass2}, $M^N$ is well-defined and a positive martingale. Moreover, $\E[Z^N_T] = 1$ for all $T \in (0,\infty)$, which implies that the positive local martingale $Z^N$ is in fact a true martingale. These facts are argued in \cite{kolliopoulos2022propagation}; see the discussion leading up to Theorem~2.4 in that paper.

%imply that $\E[ \int_0^t (\Delta B^{i,N}_s)^2 ds] < \infty$ for all $i$ and $t$, which ensures that $M^N$ is a well-defined positive martingale. We claim that $\E[Z^N_T] = 1$ for all $T \in (0,\infty)$, so that $Z^N$ is a true martingale. To see this, note that \cite[Lemma 4.3]{kolliopoulos2022propagation22} implies that for any $s < t \le T$ with $t-s$ small enough, the chaos expansion
%\[
%\frac{Z^N_t}{Z^N_s} = 1 + \sum_{m=1}^\infty \int_s^t \int_s^{t_1} \cdots \int_s^{t_{m-1}} dM^N_{t_m} \cdots dM^N_{t_1}
%\]
%converges in $L^2$. Moreover, \cite[Proposition~1]{CAKR91} together with Assumption~\ref{ass2} imply that each iterated integral has expectation zero. As a result, $\E[Z^N_t / Z^N_s] = 1$ for all such $s,t$, and this implies $\E[Z_T]=1$ as claimed.

Since $Z^N$ is a true martingale, it induces a locally equivalent probability measure $\Q^N \sim_\text{loc} \P$ under which the processes defined by
\[
W^{i,N}_t = W^i_t  - \int_0^t \Delta B^{i,N}_s ds
\]
are mutually independent standard Brownian motions. Thus under $\Q^N$ we find that $X^1,\ldots,X^N$ follow the $N$-particle dynamics \eqref{eq_particle_SDE},
\begin{equation} \label{eq_Q_particle_SDE}
dX^i_t = A(t, X^i_{[0,t]}) \left( B\left( t, X^i_{[0,t]}, \int g(t, X^i_{[0,t]}, \bm y_{[0,t]}) \mu^N_t(d\bm y) \right) dt + dW^{i,N}_t \right) + C(t, X^i_{[0,t]}) dt.
\end{equation}
We thus model the point process $\widetilde U_N$ in \eqref{pointprocess} by $U_N$~\eqref{eq_point_process_McKV} viewed under $\Q^N$, and our main result establishes convergence of the laws $\Q^N\circ U_N^{-1}$. The following is our main result.

\begin{theorem}\label{T_main_precise}
Suppose Assumptions~\ref{ass1} and \ref{ass2} are satisfied. Fix $T \in (0,\infty)$ and suppose that for some normalizing constants $a^N_T,b^N_T$ the normalized maxima of the i.i.d.\ system \eqref{eq_McKV_SDE} converge weakly to a nondegenerate distribution function $\Gamma_T$ on $\R$:
\begin{equation}\label{condition_T_main}
\P\left( \max_{i \le N} \frac{X^i_T - b^N_T}{a^N_T} \le x \right) \to \Gamma_T(x) \text{ as } N \to \infty, \quad x \in \R.
\end{equation}
Then, the sequence $(\Q^N\circ U_N^{-1})_{N \in \mathbb{N}}$ of probability laws on $M_p(E)$ converges to the law of a Poisson random measure on $(0,1]\times(x_*,x^*]$ with intensity measure $\nu$, where 
\[
\nu((a,b]\times(c,d]) = (b - a) ((1 + \gamma c)^{-1/\gamma} -(1+\gamma d)^{-1/\gamma}), \quad a<b,\ c<d.
\]
Here $x_*,x^*$ are the left and right endpoints of the support of the limiting distribution $\Gamma_T$, and $\gamma$ is the extreme value index. 
\end{theorem}

\begin{remark}\label{remark_point_process_weak_limit}
The condition~\eqref{condition_T_main} is equivalent to the weak convergence of $\P\circ U_N^{-1}$ to a Poisson random measure on on $(0,1]\times(x_*,x^*]$ with the same intensity measure $\nu$ as in Theorem~\ref{T_main_precise}. One may check this by looking at \cite[Theorem 2.1.2]{HF10}.
\end{remark}

\begin{remark}
Readers who are not familiar with the extreme value index may refer to \cite[Definition~1.1.4]{HF10}. In the case $\gamma = 0$, the expression $(1 + \gamma x)^{-1/\gamma}$ is interpreted as $\exp(-x)$. Note that the assumptions of Theorem~\ref{T_main_precise} are precisely those of the main theorem in \cite{kolliopoulos2022propagation}, even though the conclusion is stronger in that it applies to the full point process rather than to the maximum particle only.
%establish the propagation of chaos for the point process in \eqref{pointprocess}, we do not need to impose any extra condition comparing with the main theorem in \cite{kolliopoulos2022propagation22}.
\end{remark}

We now introduce some additional notation. This will prepare us to state an auxiliary result that simplifies the proof of Theorem~\ref{T_main_precise}. We write $D_3 B(t, \bm x_{[0,t]}, r)$ for the derivative with respect to $r$ and define the process
\begin{equation}\label{third_deri}
D_3 B^i_t = D_3 B\left(t, X^i_{[0,t]}, \int g(t, X^i_{[0,t]}, \bm y_{[0,t]}) \mu_t(d\bm y) \right).
\end{equation}
Note that $D_3B^i$ is adapted to the filtration $(\Fcal^{\{i\}}_t)_{t \ge 0}$ generated by $(X^i,W^i)$. We also write 
\begin{equation}\label{H_def}
H^i_t(\mu) = \int g(t, X^i_{[0,t]}, \bm y_{[0,t]}) \mu(d\bm y)
\end{equation}
for any signed measure $\mu$ on $C(\R_+)$. We then have the Taylor formula
\[
\Delta B^{i,N}_t = D_3 B^i_t \, H^i_t(\mu^N_t - \mu_t) + R^{i,N}_t \left( H^i_t(\mu^N_t - \mu_t) \right)^2,
\]
where $R^{i,N}$ is a process which is uniformly bounded in terms of the bound on the second derivative of $r \mapsto B(t, \bm x_{[0,t]}, r)$ given by Assumption~\ref{ass1}. Furthermore, we define the local martingales
\begin{equation}\label{tildem}
    \widetilde{M}_t^N = \sum_{i=1}^N \int_0^t D_3B_s^i \, H_s^i(\mu_s^N - \mu_s) dW_s^i.
\end{equation}
With all these preparations, we are ready to introduce an auxiliary theorem which is a crucial intermediate step toward establishing Theorem~\eqref{T_main_precise}. Recall that $X^1,\ldots,X^N$ are i.i.d.\ copies of the solution to \eqref{eq_McKV_SDE}, and that $Z_T^N$ is defined in \eqref{RND}. For each $m \in \N$ and $0\le s \le t$, we define the $m$-fold iterated integral $\widetilde{I}_{m}^N(s, t)$ by
\begin{equation}\label{tildei}
    \widetilde{I}_m^N(s,t) = \int_s^t\int_s^{t_1}\cdots\int_s^{t_{m-1}} d\widetilde{M}^N_{t_m}\cdots d\widetilde{M}^N_{t_1}.
\end{equation}

\begin{theorem}\label{preprocess}
Let Assumptions \ref{ass1} and \ref{ass2} be in force and fix $T \in (0,\infty)$. Let $\phi_N \colon \mathbb{R}^N\rightarrow\mathbb{R}$, $N\in\mathbb{N}$, be a sequence of measurable functions such that
\[
\sup_{N \in \N,\, (x_1,\cdots,x_N) \in \R^N}|\phi_N(x_1,\cdots,x_N)|<\infty.
\]
To show that
\begin{equation}\label{before_process}
\lim_{N\to\infty}\mathbb{E}\left[\phi_N(X_T^1,\ldots,X_T^N)(Z_T^N - 1)\right] = 0,
\end{equation}
it suffices to show that
\begin{equation}\label{T1_firstcondition}
    \lim_{N\to\infty}\left|\mathbb{E}\left[\phi_N(X_T^1,\ldots,X_T^N) \widetilde{I}_{k_1}^N(T_{0},T_1)\cdots \widetilde{I}_{k_\alpha}^N(T_{\alpha - 1},T_\alpha) \right]\right| = 0
\end{equation}
for all $\alpha \in [n]$ and all $\alpha$-tuples $(k_1,\ldots,k_\alpha)$ such that $k_\beta \in [m_\beta]\cup\{0\}$ for $\beta < \alpha$ and $k_\alpha\in[m_\alpha]$, where $m_\beta$ ($\beta=1,\ldots,\alpha)$ are some positive integers that do not depend on $N$, and $T_i = iT/n$. % are time partition on $(0,T]$ such that $T_{\alpha} - T_{\alpha - 1} = \frac{T}{n}, \alpha\in[n]$.
%\textcolor{blue}{Recall that $Z_t^N$ defined in \eqref{RND} is the Radon-Nikodym derivative define the new measure $\Q^N$. }

\end{theorem}

\section{Consequences of the main result and examples} \label{S_examples}

We now discuss some examples and applications of our main result. 

\subsection{Marginal point process}\label{marginal_point_process}
The point process $U_N$ in \eqref{eq_point_process_McKV} is defined over the product space $(0,1]\times(x_*,x^*]$. The first coordinate, which accounts for the index of the particle, is included for technical reasons: it ensures that no two points overlap. Using Theorem~\ref{T_main_precise} it is straightforward to derive a result for the marginal point process which excludes this coordinate. Indeed, by setting $I = (0,1]\times I'$ where $I'$ is a subset of $(x_*,x^*]$, we obtain that the point process defined on $(x_*,x^*]$ by
%The following is straight forward derived from our main theorem~\ref{T_main_precise}, when we set $I = [0,1]\times I'$ where $I'$ is a subset of $(x_*,x^*]$.
%If every conditions in main theorem~\ref{T_main_precise} holds, the point process defined on $(x_*,x^*]$ as:
    \[
    H_N(I') = \sum_{i = 1}^N \bm 1_{\left\{(X_T^{i} - b_T^N)/a_T^N \in I' \right\}}.
    \]
According to Theorem~\ref{T_main_precise}, the laws $\Q^N\circ H_N^{-1}$ converges weakly to a Poisson point process $H_\infty$ on $(x_*,x^*]$ with intensity measure $\mu((a,b]) = (1 + a\gamma)^{-1/\gamma} - (1 + b\gamma)^{-1/\gamma}$. The laws $\P\circ H_N^{-1}$ converges weakly to the same Poisson point process $H_\infty$ with the same intensity measure $\mu$ according to \cite[Theorem 2.1.2]{HF10}.

% We also define the point process induced by the independent copies of the solution to the limiting Mckean--Vlasov SDE \eqref{eq_McKV_SDE},
%     \[
%     \widetilde{H}_N(I') = \sum_{i = 1}^N \bm 1_{\left\{(X_T^{i} - b_T^N)/a_T^N \in I' \right\}}.
%     \]

%As a consequence of the above result, we can derive the fixed time asymptotic distribution/ joint distribution of any $k$ ranked/ top $k$ ranked particles. 

\subsection{Propagation of chaos for order statistics}
    For any $N$ particles in \eqref{eq_particle_SDE}, we use the following notation to denote the upper order statistics:
    \[
    X_T^{(1),N} \geq X_T^{(2),N} \geq X_T^{(3),N} \geq \cdots \geq X_T^{(N),N}.
    \]
As discussed in Section~\ref{S_results}, we model the $N$ populated mean--field interacting particle system~\eqref{eq_particle_SDE} by the $N$ independent copies of the solution of the Mckean--Vlasov equation~\eqref{eq_McKV_SDE} under $\Q^N$. Hence, 
we model the order statistics of the mean--field interacting particles
\[
X_T^{(1),N} \geq X_T^{(2),N} \geq X_T^{(3),N} \geq \cdots \geq X_T^{(N),N}
\]
via the order statistics of the $N$ independent copies of the solution of the Mckean--Vlasov equation
\[
X_T^{(1)} \geq X_T^{(2)} \geq X_T^{(3)} \geq \cdots \geq X_T^{(N)}
\]
under $\Q^N$. The time-$T$ marginal distribution of the rank-$k$ particle $X_T^{(k)}$ has the representation
    \[
    \mathbb{Q}^N\left(\frac{X_T^{(k)} - b_T^N}{a_T^N} \geq x\right) = \mathbb{Q}^N(H_N([x,\infty))\geq k).
    \]
    Thus, in view of Subsection~\ref{marginal_point_process} above, we have
    \[
    \lim_{N\to\infty}\Q^N\left(\frac{X_T^{(k)} - b_T^N}{a_T^N} \geq x\right) = \lim_{N\to\infty}\Q^N(H_N([x,\infty))\geq k) = \P( H_\infty([x,\infty)) \ge k),
    \]
    where $H_\infty([x,\infty))$ has a Poisson distribution with parameter $(1+x\gamma)^{-1/\gamma}$ under $\P$. The same reasoning applies to the system of independent copies, resulting in
    \[
    \mathbb{P}\left(\frac{X_T^{(k)} - b_T^N}{a_T^N} \geq x\right) = \mathbb{P}(H_N([x,\infty))\geq k),
    \]
and hence the propagation of chaos property holds for the $k$-th upper order statistic,
    \[
    \lim_{N\to\infty}\Q^N\left(\frac{X_T^{(k)} - b_T^N}{a_T^N} \geq x\right) = \lim_{N\to\infty}\mathbb{P}\left(\frac{X_T^{(k)} - b_T^N}{a_T^N} \geq x\right),
    \]
    where
    \begin{equation}\label{eq_McKean_N}
    X_T^{(1)} \geq X_T^{(2)} \geq \cdots \geq X_T^{(N)}.
    \end{equation}
    Thus, the $k$th upper order statistic from the interacting particle system \eqref{eq_particle_SDE} has the same weak limit as the $k$th upper order statistics from the independent copies of the limiting McKean--Vlasov system.
    
%    which later term is described by above example. On the other hand, for the limiting Mckean-Vlasov system \eqref{eq_McKV_SDE}, and any N copies of the system $\{X_T^i\}$, we can define the order statistics as follow:
    % \begin{equation}\label{eq_McKean_N}
    % X_T^{(1)} \geq X_T^{(2)} \geq \cdots \geq X_T^{(N)}.
    % \end{equation}
    % Since we know that $H_N(.), \tilde{H}_N(.)$ has same weak limit and
    % \[
    % \mathbb{P}(X_T^{(k)} \geq x) = \mathbb{P}(\tilde{H}_N((x,\infty))\geq k).
    % \]
    % Hence we have the propagation of chaos property for order statistics as follow:
    % \[
    % \lim_{N\to\infty}\mathbb{P}(X_T^{(k),N} \geq x) = \lim_{N\to\infty}\mathbb{P}(X_T^{(k)} \geq x)
    % \],
    % i.e the $k$-order statistics from the interacting particle system \eqref{eq_particle_SDE} has the same weak limit as the k-order statistics from the copies of the limiting McKean-Vlasov system.

\subsection{Propagation of chaos for the joint distribution of the top $k$ order statistics}\label{joint_k_order_statistics}
We now discuss the asymptotic joint distribution of the fixed time top-$k$ ranked particles of the mean--field interacting particle system,
\[
\lim_{N \to \infty} \Q^N\left( \frac{X_T^{(1)} - b^N_T}{a^N_T} \geq x_1, \ldots, \frac{X_T^{(k)} - b^N_T}{a^N_T} \geq x_k\right), \quad x_1\geq x_2 \geq \cdots \geq x_k.
\]
Note that
    \begin{align*}
        &\Q^N\left(\frac{X_T^{(1)} - b_T^N}{a_T^N} \geq x_1, \ldots, \frac{X_T^{(k)} - b_T^N}{a_T^N} \geq x_k\right)\\ 
        &= \Q^N\left(H_N([x_1,\infty))\geq 1, H_N([x_2,\infty))\geq 2, \ldots, H_N([x_k,\infty))\geq k\right)\\
        &=\sum_{\{i_1,\cdots,i_k\}\in S_k}\Q^N\left(H_N([x_1,\infty)) = i_1, H_N([x_2,x_1)) = i_2, \ldots, H_N([x_{k},x_{k-1})) = i_k  \right)
    \end{align*}
    where the set $S_k$ is defined as
    \[
    S_k = \{(i_1,\ldots,i_k): i_1 + \ldots + i_j \geq j ~~ \forall j \in [k] \}.
    \]
    Moreover, if we consider the order statistics \eqref{eq_McKean_N} of the $N$ copies of the solution to the limiting McKean--Vlasov equation, we have the following equation which establish the propagation of chaos property for the top $k$ order statistics. 
    \begin{align*}
        \lim_{N\to\infty}\mathbb\Q^N&\left(\frac{X_T^{(1)} - b_T^N}{a_T^N} \geq x_1, \ldots, \frac{X_T^{(k)} - b_T^N}{a_T^N} \geq x_k\right)\\
        &= \lim_{N\to\infty}\sum_{\{i_1,\ldots,i_k\}\in S_k }\Q^N\left(H_N([x_1,\infty)) = i_1, H_N([x_2,x_1)) = i_2, \ldots, H_N([x_{k},x_{k-1})) = i_k  \right)\\
        & = \sum_{\{i_1,\ldots,i_k\}\in S_k} \P\left(H_\infty[x_1,\infty)) = i_1, H_\infty([x_2,x_1)) = i_2 \ldots, H_\infty([x_{k},x_{k-1})) = i_k  \right)\\
        & =\lim_{N\to\infty}\P\left(\frac{X_T^{(1)} - b_T^N}{a_T^N} \geq x_1, \ldots, \frac{X_T^{(k)} - b_T^N}{a_T^N} \geq x_k\right).
    \end{align*}
    On the other hand, for each fixed tuple $(i_1,\ldots,i_k)$, as discussed in Section~\ref{marginal_point_process}, we know that the limiting distribution of $H_N$ is a Poisson point process. Thus for any non-overlapping intervals $I_1, I_2$, $H_\infty(I_1)$ is independent from $H_\infty(I_2)$. Hence we have
    \begin{align*}
    \lim_{N\to\infty}&\Q^N\left(H_N([x_1,\infty)) = i_1, H_N([x_2,x_1)) = i_2 \cdots, H_N([x_{k},x_{k-1})) = i_k  \right)\\
    &= \P\left(H_\infty([x_1,\infty)) = i_1, H_\infty([x_2,x_1)) = i_2 \cdots, H_\infty([x_{k},x_{k-1})) = i_k  \right)\\
    & = \prod_{j = 1}^k\P\left(H_\infty([x_j,x_{j-1}) = i_j\right)\\
    & = \prod_{j = 1}^k \frac{{\lambda_j}^{i_j} \e^{-\lambda_j}}{i_j!}
    \end{align*}
    where $\lambda_j = (1 + \gamma x_j)^{-1/\gamma} - (1 + \gamma x_{j - 1})^{-1/\gamma}$.
Thus we have:
    \begin{equation} \label{eq_upper_k_order_statistics_limit}
    \lim_{N\to\infty}\Q^N\left(\frac{X_T^{(1)} - b_T^N}{a_T^N} \geq x_1, \cdots, \frac{X_T^{(k)} - b_T^N}{a_T^N} \geq x_k\right) = 
    \sum_{\{i_1,\ldots,i_j\}\in S_k}\prod_{j = 1}^k \frac{{\lambda_j}^{i_j} \e^{-\lambda_j}}{i_j!}.
    \end{equation}
There is another way to characterize this limit. Due to \cite[Theorem 2.1.1]{HF10}, the normalized top $k$ ranked particles of the independent system,
    \[
    \left(\frac{X_T^{(1)} - b_T^N}{a_T^N}, \frac{X_T^{(2)} - b_T^N}{a_T^N}, \ldots, \frac{X_T^{(k)} - b_T^N}{a_T^N} \right),
    \]
    converge in distribution to
    \[
    \left(\frac{(E_1)^{-\gamma} - 1}{\gamma}, \frac{(E_1 + E_2)^{-\gamma} - 1}{\gamma}, \cdots, \frac{(E_1 + E_2 + \cdots + E_k)^{-\gamma} - 1}{\gamma}\right),
    \]
    where $\gamma$ is the extreme value index, and $E_i$ are i.i.d exponential distribution.
Theorem~\ref{T_main_precise} implies that this limit coincides with the one in \eqref{eq_upper_k_order_statistics_limit}.

% There is another characterization of the joint distribution of top $k$-ranked particles, see \cite[Theorem 2.1.1]{HF10}. The only thing left us to show here is the propagation of chaos property i.e the limiting joint distribution of top $k$ ranked particles of interacting system is the same as the independent one. The joint distribution of the re-normalzied top $k$ ranked particles of the independent system:
%     \[
%     \left(\frac{X_T^1 - b_T^N}{a_T^N}, \frac{X_T^2 - b_T^N}{a_T^N}, \cdots, \frac{X_T^k - b_T^N}{a_T^N} \right),
%     \]
%     converge in distribution to:
%     \[
%     \left(\frac{(E_1)^{-\gamma} - 1}{\gamma}, \frac{(E_1 + E_2)^{-\gamma} - 1}{\gamma}, \cdots, \frac{(E_1 + E_2 + \cdots + E_k)^{-\gamma} - 1}{\gamma}\right),
%     \]
%     where $\gamma$ is the extreme distribution parameter, and $E_i$ are i.i.d exponential distribution.

\subsection{Propagation of chaos when $k\to\infty$}

Let us comment briefly on the \emph{intermediate order statistics} of the interacting particle system \eqref{eq_particle_SDE}. The intermediate order statistics are particles $X^{k(N),N}_t$ where the rank $k(N)$ increases with $N$ in such a way that $k(N)/N \to 0$ as $N\to \infty$. Let us consider the general case where $k:= k(N)\to\infty$. Our previous computations show that
%Following the spirits in the example of order statistics, we are now discussing the intermediate order statistics of the interacting particles system \eqref{eq_particle_SDE}. The intermediate order statistics is nothing else but let the parameter $k$ be depend on number of the particles N, and has relation $k(N)/N \to \rho,~N\to \infty$, where $\rho\in(0,1)$. Let's consider a simple case where $k(N) = \rho N$. Note we have the following equation:
\begin{align*}
\lim_{N\to\infty}\Q^N\left(\frac{X_T^{k(N)} - b_T^N}{a_T^N} \geq x\right) &= \lim_{N\to\infty}\Q^N\left(H_N([x,\infty))\geq k(N)\right)\\
    &= \lim_{N\to\infty}\Q^N\left(\frac{H_N([x,\infty))}{k(N)}\geq 1\right)\\
    &=0.
\end{align*}
This indicates that our results do not yield any information about the asymptotic behavior of the $k(N)$-th order statistics. In particular, our results do not yield any information about the asymptotic behavior of intermediate order statistics as well. In fact, it is well-known in the classical extreme value literature that certain second-order conditions are required to establish convergence of suitably normalized intermediate order statistics. An interesting open question is whether those conditions are also sufficient to obtain the analogous result for the interacting systems that we consider here, but this is outside the scope of the current paper.

%This suggests that under the origin assumptions, we are not able to deduced the propagation of chaos property for the intermediate order statistics. When one look back at the literature about the asymptotic behavior of intermediate order statistics (e.g in \cite[Chapter 2]{HF10}), except the condition that there exists re-normalizing factor $a_N, b_N$, there are also second order conditions which allow us to establish the convergence of the intermediate order statistics for i.i.d random variable. Naturally, one may expect we need such condition to establish the propagation of chaos property for the intermediate order statistics which beyond the discussion of this paper.

\subsection{Two examples}

We now discuss two examples that illustrate the main result as well as it's applications.

\begin{example}[Gaussian particles] \label{ex_1}
We build on \cite[Example 3.1]{kolliopoulos2022propagation}. Consider the following Gaussian particle system:
\begin{align*}
X_{t}^{i, N} &= X_{0}^{i} - \kappa \int_{0}^{t}\left(X_{s}^{i, N} -  \frac{1}{N}\sum_{j=1}^NX_{s}^{j, N} \right)ds + \sigma W_t^{i}, \quad i=1,\ldots,N,
\end{align*}
with i.i.d.\ $N(m_0,\sigma_0^2)$ initial conditions. Here $\kappa, m_0 \in \R$ and $\sigma, \sigma_0 \in (0,\infty)$ are parameters. In our setting this example arises by taking $A(t,\bm x_{[0,t]}) = \sigma$, $B(t,\bm x_{[0,t]},r) = -\kappa (x_t - r)/\sigma$, $C(t, \bm x_{[0,t]}) = 0$, and $g(t,\bm x_{[0,t]},\bm y_{[0,t]}) = y_t$. Clearly Assumption~\ref{ass1} is satisfied. The McKean--Vlasov equation \eqref{eq_McKV_SDE} reduces to
\[
dX_t = - \kappa \left(X_t - \E[X_t] \right)dt + \sigma dW_t.
\]
From the discussion in \cite[Example 3.1]{kolliopoulos2022propagation} we know that for suitable parameters $a^N_T, b^N_T$,
\[
\P\left( \max_{i \le N} \frac{X^i_T - b^N_T}{a^N_T} \le x \right)
%= \P\left( \max_{i \le N} \frac{(X^i_T - m_0)/\sigma_T - b^N}{a^N} \le x \right)
\to \exp(-e^{-x}) \text{ as } N \to \infty.
\]
That is, $X_T$ belongs to the Gumbel domain of attraction. Our main result Theorem \ref{T_main_precise} applies, and by the discussion in Subsection~\ref{marginal_point_process}, we know that the marginal point process converges weakly to a Poisson random measure on the real line $\mathbb{R}$ with intensity measure $\mu((a,b]) = \exp(-a) - \exp(-b)$. Thanks to \eqref{eq_upper_k_order_statistics_limit} we have the following asymptotics:
    \[
    \lim_{N\to\infty}\P\left(X_T^{(1),N} \geq x_1, \ldots, X_T^{(k),N} \geq x_{k}\right) = 
    \sum_{\{i_1,\ldots,i_j\}\in S_k}\prod_{j = 1}^k \frac{{\lambda_j}^{i_j} \e^{-\lambda_j}}{i_j!},
    \]
where $\lambda_j = \exp(-x_j) - \exp(-x_{j-1})$ and the set $S_k$ is defined as
    \[
    S_k = \{(i_1,\ldots,i_k): i_1 + \cdots + i_j \geq j ~~ \forall j \in [k] \}.
    \]
\end{example}

\begin{example}[Rank-based diffusions] \label{ex_2}
This is another example from \cite{kolliopoulos2022propagation}. Consider the $N$-particle system evolving according to
\[
dX^{i,N}_t = B\left( F^N_t( X^{i,N}_t) \right) dt + \sqrt{2} dW^i_t
\]
for $i \in [N]$, where $B(r)$ is a twice continuously differentiable function on $[0,1]$ and
\[
F^N_t(x) = \frac{1}{N} \sum_{j=1}^N \bm1_{\{X^{j,N}_t \le x\}}
\]
is the empirical distribution function. If X is stationary, an sufficient condition for the stationary distribution of the limiting Mckean--Vlasov equation was given in Example~3.2 in \cite{kolliopoulos2022propagation}. Under this condition, the limiting stationary distribution belongs to the Gumbel domain of attraction, which is exactly the same as in the previous example. Hence we may apply the main Theorem~\ref{T_main_precise} to derive the weak limit of the marginal point process as well as the asymptotic joint distribution of the top $k$ order statistics.
\end{example}

\section{Proof of Theorem~\ref{preprocess}}
\label{preprocess_proof}

Theorem~\ref{preprocess} is a technical result that encapsulates a significant portion of the technical work required to prove Theorem~\ref{T_main_precise}. The crucial observation is that parts of the proof of \cite[Theorem~2.4]{kolliopoulos2022propagation} go through in a substantially generalized form, which we make use of here. In order to avoid repeating several pages of technical computations and estimates from \cite{kolliopoulos2022propagation} essentially verbatim, we simply describe the minor adjustments required to establish Theorem~\ref{preprocess}. To assist the reader in following the argument, we use the same notation as in \cite{kolliopoulos2022propagation}.

To prove Theorem~\ref{preprocess} we repeat, word for word, the argument in Step~1--3 and the first part of Step 4, leading up to (6.13), in \cite[Section~6]{kolliopoulos2022propagation}, with one change:
\[
\prod_{i=1}^N \bm 1_{\{X_T^i < x_N\}}
\]
in that proof is replaced with 
\[
\phi_N(X_T^1,\cdots,X_T^N).
\]
This establishes that there exists $n$ such that for any $\epsilon > 0$, we have
\[
\E\left[\phi_N(X_T^1,\ldots,X_T^N)(Z_T^N - 1)\right] = \sum_{\alpha = 1}^n A_\alpha^N,
\]
where for some constant $C$,
\begin{equation}\label{eq_bound_A_alpha}
\left|A_\alpha^N\right|\leq \alpha\epsilon + \frac{C}{\sqrt{N}} + \sum_{k_1,\ldots,k_\alpha} \left|\mathbb{E}\left[\phi_N(X_T^1,\ldots,X_T^N)\widetilde{I}_{k_1}^N(T_{0},T_1)\cdots \widetilde{I}_{k_\alpha}^N(T_{\alpha - 1},T_\alpha) \right]\right|.
\end{equation}
Here $\widetilde{I}_{m}^N(s,t)$ is defined in \eqref{tildei} and $T_i = iT/n$. The sum in \eqref{eq_bound_A_alpha} ranges over all $(k_1,\ldots, k_\alpha)$ such that $k_\beta \in [m_{\beta}]\cup \{0\}$ for $\beta < \alpha$ and $k_{\alpha}\in[m_\alpha]$, where $m_\beta$ ($\beta=1,\ldots,\alpha)$ are positive integers that do not depend on $N$.

Now, note that $\epsilon > 0$ is arbitrary and that $n$, $C$, and the number of terms in the sum in \eqref{eq_bound_A_alpha} do not depend on $N$. Hence to show that \eqref{before_process} holds, it suffices to show that \eqref{T1_firstcondition} holds for all $\alpha \in [n]$ and all $\alpha$-tuples $(k_1,\ldots,k_\alpha)$ such that $k_\beta \in [m_\beta]\cup\{0\}$ for $\beta < \alpha$ and $k_\alpha\in[m_\alpha]$. This completes the proof.

\section{Proof of Theorem~\ref{T_main_precise}} \label{S_main_proof}

To prove Theorem~\ref{T_main_precise}, we employ a variation of \cite[Proposition 3.22]{MR2364939} which provides conditions under which a sequence $\{\xi_n\}_{n \in \mathbb{N}}$ of point processes converges weakly to a candidate limit $\xi$. Observe that the proof of that result uses \cite[Equation (3.24)]{MR2364939} only for showing that the sequence $\{\xi_n\}_{n \in \mathbb{N}}$ is tight and for obtaining the inequality
\begin{equation}\label{eq_expectation_inequality}
\E[\eta(I)] \leq \E[\xi(I)]
\end{equation}
for any weak limit $\eta$ of $\{\xi_n\}_{n \in \mathbb{N}}$ (through some subsequence) and any set $I$ taken from a certain collection $\Tcal$ of subsets of the state space $E$. The inequality \eqref{eq_expectation_inequality} for all $I \in \Tcal$ is established at the end of the proof as part of showing that the arbitrary weak limit $\eta$ is a simple point process. Therefore, we are going to replace \cite[Equation (3.24)]{MR2364939} with another condition that implies both the tightness of $\{\xi_n\}_{n \in \mathbb{N}}$ and \eqref{eq_expectation_inequality}. 

We start by noting that we only need \eqref{eq_expectation_inequality} to hold as an equality, which is equivalent to
\begin{equation*}
\P(\eta(I) = t) = \P(\xi(I) = t)
\end{equation*} 
for each $t \in \mathbb{N}$ and all $I \in \Tcal$. By the continuity argument immediately after \cite[Display~(3.28)]{MR2364939}, the latter reduces to having 
\begin{equation}\label{limeq_prob_equality}
\lim_{n \to \infty}\P(\xi_n(I) = t) = \P(\xi(I) = t)
\end{equation}
for all $I \in \Tcal$ and $t \in \mathbb{N}$. We claim now that \eqref{limeq_prob_equality} is also sufficient for tightness. Indeed, by \cite[Exercise 3.5.2]{MR2364939}, we only need to show that
\begin{equation}\label{tightness}
\lim_{t\to\infty}\limsup_{n\to\infty}\P(\xi_n(K) > t) = 0
\end{equation}
for any compact subset $K$ of the state space $E$. Thanks to the properties of $\Tcal$, one can pick a set $I \in \Tcal$ such that $K\subset I$, so we can bound
\begin{align*}
\limsup_{N\to\infty}\P(\xi_n(K) > t) &\leq \limsup_{N\to\infty}\P(\xi_n(I) > t) \\
    & = \P(\xi(I) > t),
\end{align*}
where the last equality follows from \eqref{limeq_prob_equality}. Sending $t\to\infty$ in the above estimate we get \eqref{tightness} so the required tightness follows. 

Observe now that \cite[Equation (3.23)]{MR2364939} is also obtained from  \eqref{limeq_prob_equality} by simply setting $t = 0$.
Therefore, one can deduce the weak convergence $\xi_n \to \xi$ as $n \to \infty$ by verifying the following two conditions: 
\begin{enumerate}[label=\textnormal{(\alph*)}]
    \item $\P(\xi(\partial I) = 0) = 1$ for all $I \in \Tcal$. %which is true in our case since $\xi = U_\infty$ is a Poisson point process.
    \item $\displaystyle{\lim_{n\to\infty}\P(\xi_n(I) = t ) = \P(\xi(I) = t)}$ for all $I\in \Tcal$ and $t \in \mathbb{N}$.
\end{enumerate}
In our case, $\xi = U_\infty$ is a Poisson point process with intensity $\nu$ so that $\P \circ U_{N}^{-1} \to \text{Law}(U_\infty)$ as $N \to \infty$, $\{\xi_n\}_{n \in \mathbb{N}}$ is a sequence of point processes $\{\widetilde U_N\}_{N \in \mathbb{N}}$ with $\Q^N\circ U_N^{-1}$ being the law of $\widetilde U_N$ for all $N \in \mathbb{N}$, $E = (0,1]\times(x_*, x^*]$, and finally $\Tcal = \Ical$ where $\Ical$ is the family of finite unions of rectangles of the form:
\[
I = \bigcup_{j=1}^m(a_j,b_j]\times(c_j,d_j],
\]
where $m \in \N$, $0 < a_j < b_j \leq 1$ and $x_* < c_j < d_j \leq x^*$ for $j \in [m]$. Without loss of generality, we assume that all the rectangles that constitute each $I$ are pairwise disjoint.
Since $\xi = U_\infty$ is a Poisson point process with intensity measure that is absolutely continuous with respect to the Lebesgue measure, (a) follows immediately. 
Then, since $\P \circ U_{N}^{-1} \to \text{Law}(U_\infty)$ as $N \to \infty$, condition (b) translates to
\begin{equation}\label{T2_claim_1'}
\lim_{N\to\infty}\Q^N({U}_N(I) = t) = \lim_{N\to\infty}\P({U}_N(I) = t),
\end{equation} 
for any set $I \in \Ical$ and any $t \in \mathbb{N}$.
Therefore, the proof of our result reduces to showing that \eqref{T2_claim_1'} is true for any set $I \in \Ical$ and any $t \in \mathbb{N}$. This will be done by employing ideas from \cite{kolliopoulos2022propagation}.

\paragraph{Proof of condition~\eqref{T2_claim_1'}:}
Observe that what we want to show now is equivalent to
\begin{equation}\label{T2_claim_result}
    \lim_{N\to\infty} \mathbb{E}\left[\bm 1_{\{{U}_N(I) = t\}}(Z_T^N - 1)\right] = 0,
\end{equation}
for any $I \in \Ical$ and any integer $t \geq 0$. Recalling Theorem~\ref{preprocess} with $\phi_N(X_T^1,\ldots,X_T^N)$ taken to be the indicator $\bm 1_{\{{U}_N(I) = t\}}$ which is a bounded and measurable function of $X_T^1,\ldots,X_T^N$, it suffices to show that 
\begin{equation}\label{T2_claim2}
\lim_{N\to\infty}\left|\mathbb{E}\left[\bm 1_{\{{U}_N(I) = t\}} \widetilde{I}_{k_1}^N(T_{0},T_1)\cdots \widetilde{I}_{k_\alpha}^N(T_{\alpha - 1},T_\alpha) \right]\right| = 0.
\end{equation}
Under the convention $$I_N = \bigcup_{j=1}^m(Na_j, Nb_j] \times (a_T^Nc_j + b_T^N, a_T^Nd_j +b_T^N]$$ and also $$\mathbf{X}_T^i = (i, X_T^i),$$
%\[
%K_N:= \left\{i:\frac{i}{N} \in (a, b]\right\} %\subset[N],
%\]
%\[
%A_N:=(a_T^Nc + b_T^N, a_T^Nd +b_T^N],
%\]
we can re-write the indicator $\bm 1_{\{{U}_N(I) = t\}}$ as follows:
% \begin{align}\label{indicator-rewrite}
% \bm 1_{\{{U}_N(I) = t\}} = & \sum_{\{i_1,\cdots,i_t\}\subset[N]}\prod_{l=1}^t\bm 1_{\{X_T^{i_l} > x_N\}}\prod_{i\in[n]/\{i_1,\ldots,i_t\}} \bm 1_{\{X_t^i \leq x_N\}}\\
% & = \sum_{\{i_1,\cdots,i_t\}\subset[N]}  \sum_{\kappa = 0}^{N-t}(-1)^\kappa\prod_{l = 1}^t \bm 1_{\{X_T^{i_l} > x_N\}} \sum_{{\{i'_1,\cdots,i'_\kappa\}\subset[N]/\{i_1,\ldots,i_t\}}}\prod_{l=1}^\kappa \bm 1_{\{X_T^{i'_l} > x_N\}}.
% \end{align}
\begin{align}\label{indicator-rewrite}
\bm 1_{\{{U}_N(I) = t\}} = & \sum_{\{i_1,\cdots,i_t\}\subset [N]}\prod_{l=1}^t\bm 1_{\{\mathbf{X}_T^{i_l} \in I_N\}}\prod_{i\in [N]/\{i_1,\ldots,i_t\}} \bm 1_{\{\mathbf{X}_T^i \notin I_N\}}\nonumber \\
& = \sum_{\{i_1,\cdots,i_t\}\subset [N]} \prod_{l = 1}^t \bm 1_{\{\mathbf{X}_T^{i_l} \in I_N\}} \sum_{\kappa = 0}^{N-t}(-1)^\kappa \sum_{{\{i'_1,\cdots,i'_\kappa\}\subset [N]/\{i_1,\ldots,i_t\}}}\prod_{l=1}^\kappa \bm 1_{\{\mathbf{X}_T^{i'_l} \in I_N\}} \nonumber \\
& = \sum_{\{i_1,\cdots,i_t\}\subset [N]} \sum_{\kappa = 0}^{N-t}(-1)^\kappa \sum_{{\{i'_1,\cdots,i'_\kappa\}\subset [N]/\{i_1,\ldots,i_t\}}}\prod_{l = 1}^t \bm 1_{\{\mathbf{X}_T^{i_l} \in I_N\}} \prod_{l=1}^\kappa \bm 1_{\{\mathbf{X}_T^{i'_l} \in I_N\}},
\end{align}
and we proceed now to the expansion of the iterated integrals $\widetilde{I}_{k_{\beta}}^N(T_{\beta - 1}, T_{\beta})$ in equation \eqref{T2_claim2}. In view of \eqref{tildem} and the definitions of $H_s^i(\mu^N_s - \mu_s)$ and $\mu^N_s$ in \eqref{H_def} and \eqref{iem} respectively, we can write
\begin{align}\label{exhaustiveexpansion}
\widetilde{M}_t^N = \frac{1}{N}\sum_{i=1}^N\sum_{j=1}^N\int_0^tG_s^{ij} dW_s^i,
\end{align}
where
\begin{equation}\label{G}
    G_t^{ij} = D_3B_t^i\, \left(g\left(t,X_{[0,t]}^i,X_{[0,t]}^j\right) - \int g\left(t,X_{[0,t]}^i,\bm y_{[0,t]}\right) \mu_t(d{\bm y}) \right)
\end{equation}
for all $i, j \in [N]$. Recall that $D_3B_t^i$ is the partial derivative of drift function $B(t, \bm x_{[0,t]}, r)$ with respect to its third argument as defined in \eqref{third_deri}. Plugging \eqref{exhaustiveexpansion} and \eqref{G} into the definition of $\widetilde{I}_{k_{\beta}}^N(T_{\beta - 1}, T_{\beta})$, see \eqref{tildei}, gives
\begin{equation}\label{tilde_I}
    \widetilde{I}^N_{k_\beta}(T_{\beta - 1}, T_\beta) = \frac{1}{N^{k_\beta}}\sum_{{\bm i}\in[N]^{k_\beta}}\sum_{{\bm j}\in[N]^{k_\beta}}I_{{\bm i},{\bm j}}^{N}(T_{\beta - 1}, T_\beta),
\end{equation}
where we use the iterated integral notation from \cite[equation (4.1)]{kolliopoulos2022propagation}, namely:

\begin{equation} \label{iterated_integral_I}
    I_{{\bm i},{\bm j}}^{N}(s, t) = \int_s^t G_{t_1}^{i_1,j_1} \int_s^{t_1} G_{t_2}^{i_2,j_2}\cdots\int_s^{t_{k-1}}G_{t_k}^{i_k,j_k} dW_{t_k}^{i_k}\cdots dW_{t_1}^{i_1}
\end{equation}
for any $k \in \N$ and any multiindices ${\bm i} = (\ir_1,\ldots,\ir_k) \in [N]^k$ and ${\bm j} = (\jr_1,\ldots,\jr_k) \in [N]^k$. ~Using the representation of $\widetilde{I}^N_{k_\beta}(T_{\beta - 1}, T_\beta)$ in \eqref{tilde_I}, the product of iterated integrals appearing in \eqref{T2_claim2} can be written as
\begin{equation}\label{T2_prod_iterated_integral}
\widetilde{I}_{k_1}^N(T_{0},T_1)\cdots \widetilde{I}_{k_\alpha}^N(T_{\alpha - 1},T_\alpha)
= \frac{1}{N^{k_1 + \cdots + k_\alpha}}\sum_{({\bm i}_1,{\bm j}_1),\ldots,({\bm i}_\alpha,{\bm j}_\alpha)} \prod_{\beta = 1}^{\alpha} I_{{\bm i_{\beta}},{\bm j_{\beta}}}^{N}(T_{\beta - 1}, T_\beta).
\end{equation}
Now we can combine the identities \eqref{T2_prod_iterated_integral} and \eqref{indicator-rewrite} and apply the triangle inequality to get
\begin{align}\label{T2_indicator_inequality}
\Big|\mathbb{E}&\left[\bm 1_{\{{U}_N(I) = t\}}\widetilde{I}_{k_1}^N(T_{0},T_1)\cdots \widetilde{I}_{k_\alpha}^N(T_{\alpha - 1},T_\alpha) \right]\Big| \nonumber\\
&\leq \frac{1}{N^s}\sum_{\{i_1,\ldots,i_t\}\subset [N]}\sum_{\kappa=0}^{N-t}\sum_{\{i'_1,\ldots,i'_\kappa\}\subset [N]/\{i_1,\ldots,i_t\}}\sum_{({\bm i}_1,{\bm j}_1),\ldots,({\bm i}_\alpha,{\bm j}_\alpha)}\nonumber\\
&\qquad \qquad \qquad \qquad \left|\mathbb{E}\left[\prod_{l=1}^t \bm 1_{\{\mathbf{X}_T^{i_l} \in I_N\}}\prod_{l=1}^\kappa\bm 1_{\{\mathbf{X}_T^{i'_l} \in I_N\}} \prod_{\beta = 1}^{\alpha} I_{{\bm i_{\beta}},{\bm j_{\beta}}}^{N}(T_{\beta - 1}, T_\beta)\right]\right|\nonumber \\
&= \frac{1}{N^s}\sum_{\kappa=0}^{N-t}\sum_{\{i_1,\ldots,i_t\}\subset[N]}\sum_{\{i'_1,\ldots,i'_\kappa\}\subset[N]/\{i_1,\ldots,i_t\}}\sum_{({\bm i}_1,{\bm j}_1),\ldots,({\bm i}_\alpha,{\bm j}_\alpha)}\nonumber\\
&\qquad \qquad \qquad \qquad\left|\mathbb{E}\left[\prod_{l=1}^t \bm 1_{\{\mathbf{X}_T^{i_l} \in I_N\}}\prod_{l=1}^\kappa\bm 1_{\{\mathbf{X}_T^{i'_l} \in I_N\}} \prod_{\beta = 1}^{\alpha} I_{{\bm i_{\beta}},{\bm j_{\beta}}}^{N}(T_{\beta - 1}, T_\beta)\right]\right|,
\end{align}
where $s = k_1 + \cdots + k_\alpha$. Therefore, to show \eqref{T2_claim2}, it suffices to show the last sum in the above goes to zero as $N\to\infty$.

%We make a pause here to see what we have so far. We have established a sufficient condition via \eqref{T2_indicator_inequality} to conclude the tightness of the law of $\Q^N\circ{U}_N^{-1}$. To estimate the right hand side of the \eqref{T2_indicator_inequality}, we will apply \cite[Lemma 4.1]{kolliopoulos2022propagation22} to get an upper bound for the number of non-zero terms appearing in the summation of right hand side of \eqref{T2_indicator_inequality}. For all those non-zero terms, we will derive a universal bound. Combing these upper bounds, we will get an upper bound of right hand side of \eqref{T2_indicator_inequality} which goes to zero as $N\to\infty$.

We count now the number of non-zero terms in the last sum in \eqref{T2_indicator_inequality}. Note that each term in that sum has the form 
\begin{equation}\label{crucialterm}
\E\left[ \Psi \prod_{\beta = 1}^{\alpha} I_{{\bm i_{\beta}},{\bm j_{\beta}}}^{N}(T_{\beta - 1} , T_{\beta}) \right]
\end{equation}
where 
% \[
% \Psi = \prod_{l=1}^t \bm 1_{\{X_T^{i_l} > x_N\}}\prod_{l=1}^\kappa \bm 1_{\{X_T^{i'_l} > x_N\}},
% \]
\[
\Psi = \prod_{l=1}^t \bm 1_{\{X_T^{i_l} \in A_N\}}\prod_{l=1}^\kappa \bm 1_{\{X_T^{i'_l} \in A_N\}},
\]
is a bounded $\Fcal_T^K$-measurable random variable for $K:=\{i_1,\ldots,i_t,i'_1,\ldots,i'_\kappa\}$, each iterated integral $I_{{\bm i_{\beta}},{\bm j_{\beta}}}^{N}(T_{\beta - 1} , T_{\beta})$ is given by equation \eqref{iterated_integral_I} in terms of the processes $\nG^{\ir\jr}$, and we can verify that for all $V \subset [N]$, $\ir \in V$, $\jr \notin V$ we have:
\begin{align*}
\E\left[ \nG^{\ir\jr}_s \mid \Fcal_t^V \right] &= \E \left[ D_3B_t^i\,\left(g\left(t,X_{[0,t]}^i,X_{[0,t]}^j\right) - \int g\left(t,X_{[0,t]}^i,\bm y_{[0,t]}\right) \mu_t(d{\bm y}) \right)\mid \Fcal_t^V \right]\\
&=D_3B_t^i\left(\E\left[g\left(t,X_{[0,t]}^i,X_{[0,t]}^j\right) \mid\Fcal_t^V \right] - \E \left[g\left(t,X_{[0,t]}^i,X_{[0,t]}^j\right) \mid\Fcal_t^V \right]\right) \\
&= 0
\end{align*}
Thus, to give an upper bound for the number of terms that survive in the last sum in  \eqref{T2_indicator_inequality}, we need to count the terms in that sum for which $K$ violates both Conditions 1 and 2 of ~\cite[Lemma 4.1]{kolliopoulos2022propagation}. By \cite[Lemma 4.2]{kolliopoulos2022propagation}, there are at most 
\[
\binom{N}{\kappa + t} (\kappa + t) (\kappa + t + 1) \cdots (\kappa + t + s) N^{s-1},
\]
choices of $K$ with $|K| = \kappa + t$ for which \eqref{crucialterm} is not zero. Moreover, the last sum in \eqref{T2_indicator_inequality} consists of terms of the form \eqref{crucialterm} for $|K| \geq t$, where each term with $|K| = \kappa + t$ for $\kappa \geq 0$ appears precisely
\[
\binom{t + \kappa}{t} = \binom{t + \kappa}{\kappa}
\]
times. The latter follows from the fact that the term in \eqref{crucialterm} for that particular $K$ emerges exactly once for each choice of the $\kappa$ indices $i'_1,\ldots,i'_\kappa$ for which $\{i'_1,\ldots,i'_\kappa\} \subset K$, by picking then the remaining $t$ indices $\{i_1,\ldots,i_t\}$. Hence, there are at most 
\begin{equation}\label{number_of_terms_survive}
\binom{t + \kappa}{t}\binom{N}{\kappa + t} (\kappa + t) (\kappa + t + 1) \cdots (\kappa + t + s) N^{s-1},
\end{equation}
surviving terms of the form \eqref{crucialterm} in the last sum in \eqref{T2_indicator_inequality} with $|K| = t + \kappa$

% Now we are deriving a universal bound for each survival term regardless of the choice of fixed set $\{i_1,\ldots,i_t\}$ in \eqref{T2_indicator_inequality}. Under the convention $x_N = a_T^Nx + b_T^N$. From the assumption we have 
% \[
% \left(1 - \P(X_T^i > x_N)\right)^N = \P\left(\max_{i = 1,\cdots,N} X_T^i < x_N\right) \to \Gamma_T(x) > 0,
% \]
% which implies 
% \[
% N \P(X_T > x_N) \le - N \log\left(1-\P(X_T > x_N)\right) \le C.
% \] This gives us $L^p$ bound for the indicator function part in \eqref{T2_indicator_inequality}
% \[
% P(X_T^i > x_N)^p \leq (C/N)^p.
% \]
Next, we are going to derive a common bound for all surviving terms in \eqref{T2_indicator_inequality}. Note that for any $(i, x) \in I_N$, one has $x > a_T^Nc + b_T^N$ where $c = \min\{c_j: 1 \leq j \leq m\}$. Moreover, from our assumptions we have 
\begin{align*}
\left(1 - \P(X_T^1 > a_T^Nc + b_T^N)\right)^N &= \P\left(\max_{i \leq N} X_T^i < a_T^Nc + b_T^N\right) \\
&= \P\left(\frac{\max_{i \leq N}X_T^i - b_T^N}{a_T^N} < c\right) \\
&\to \Gamma_T(c) > 0
\end{align*}
as $N \to \infty$, which implies 
\[
N \P(X_T^1 > a_T^Nc + b_T^N) \le - N \log\left(1-\P(X_T^1 > a_T^Nc + b_T^N)\right) \le C
\] 
for all $N$ and some $C > 0$. Hence, we can combine the above to obtain
\begin{align}\label{probound}
\P(\mathbf{X}_T^i \in I_N) &\leq P(X_T^i > a_T^Nc + b_T^N) \nonumber \\
&= P(X_T^1 > a_T^Nc + b_T^N)\nonumber \\
&\leq \frac{C}{N}.
\end{align}
Furthermore, the estimations in \cite[Page 30]{kolliopoulos2022propagation} give
\begin{equation}\label{iterated_integral_estimation_tightness}
\left\|  I_{{\bm i_{\beta}},{\bm j_{\beta}}}^{N}(T_{\beta - 1}, T_\beta)  \right\|_p\leq C (4\sqrt{p})^{k_\beta} \prod_{\ell=1}^{k_\beta} \sqrt{p} = C p^{k_\beta}.
\end{equation}
for some arbitrary constant $C$. Therefore, we can apply H\"older's inequality with some exponents $p_N$ and $q_N = p_N/(p_N - 1)$ and use also \eqref{probound} to obtain
\begin{align}\label{universal_bound}
\Bigg\vert \mathbb{E} &\left[ \prod_{l=1}^t \bm 1_{\{\mathbf{X}_T^{i_l} \in I_N\}}\prod_{l=1}^\kappa\bm 1_{\{\mathbf{X}_T^{i'_l} \in I_N\}} \prod_{\beta = 1}^{\alpha} I_{{\bm i_{\beta}},{\bm j_{\beta}}}^{N}(T_{\beta - 1}, T_\beta) \right] \Bigg\vert \nonumber \\
&\le
\mathbb{E} \left[\prod_{l=1}^t \bm 1_{\{\mathbf{X}_T^{i_l} \in I_N\}}\prod_{l=1}^\kappa\bm 1_{\{\mathbf{X}_T^{i'_l} \in I_N\}} \right]^{\frac{1}{q_N}} \prod_{\beta = 1}^{\alpha} \left\| I_{{\bm i_{\beta}},{\bm j_{\beta}}}^{N}(T_{\beta - 1}, T_\beta) \right\|_{\alpha p_N}  \nonumber \\
&= \left(\prod_{l=1}^t \P\left(\mathbf{X}_T^{i_l} \in I_N\right)\prod_{l=1}^\kappa \P\left(\mathbf{X}_T^{i'_l} \in I_N\right) \right)^{\frac{1}{q_N}} \prod_{\beta = 1}^{\alpha} \left\| I_{{\bm i_{\beta}},{\bm j_{\beta}}}^{N}(T_{\beta - 1}, T_\beta) \right\|_{\alpha p_N}  \nonumber \\
&\leq C\left(\frac{C}{N}\right)^{\frac{\kappa + t}{q_N}}p_N^s.
\end{align}
which is the required common bound.

Finally, we apply on the last sum in \eqref{T2_indicator_inequality} the upper bounds \eqref{number_of_terms_survive} and \eqref{universal_bound} we have derived for the numbers and the sizes of non-zero terms respectively, so taking also $p_N = \ln N$ we can bound that sum by
\begin{align}\label{eq_final_bound}
    \frac{1}{N^s}&\sum_{\kappa = 0}^{N - t}  \binom{t + \kappa}{t}\binom{N}{\kappa + t}N^{s-1}(\kappa+t)\cdots(\kappa+t+s) C\lceil \ln N \rceil^s\left(\frac{C}{N}\right)^{(\kappa + t)(1 - \frac{1}{
    \ln N})}\nonumber\\
    & = \sum_{\kappa = t}^N \frac{1}{N}\binom{N}{\kappa}\binom{\kappa}{t}\kappa \cdots(\kappa + s)C\lceil \ln N \rceil^s\left(\frac{C}{N}\right)^{\kappa(1 - \frac{1}{\ln N})}\nonumber\\
    &\leq \frac{C\lceil \ln N \rceil^s}{N}\sum_{\kappa = 1}^N \binom{N}{\kappa}\kappa \cdots (\kappa + s + t)\left(\frac{C}{N}\right)^{\kappa(1 - \frac{1}{\ln N})}\nonumber \\
    &\leq \frac{C\lceil \ln N \rceil^s}{N}(s+t+2)(s+t+1)^{2(s+t+1)}e^{Ce}(ce)^{s+t+1}
\end{align}
In the above display, the last inequality came from applying \cite[Lemma 4.6]{kolliopoulos2022propagation} with $S = s + t$, and we have also used the bound
% identity
% \[
% \binom{N}{t}\binom{N-t}{\kappa - t} = \binom{N}{\kappa}\binom{\kappa}{t}
% \]
% to get the first equality, and the 
\[
\binom{\kappa}{t}\kappa\cdots(\kappa+s) \leq \kappa\cdots(\kappa+s+t).
\]
Since $s = k_1 + \cdots + k_\alpha$ and $t$ do not depend on $N$, the right hand side of \eqref{eq_final_bound} vanishes as $N\to\infty$, which means that the last sum in \eqref{T2_indicator_inequality} goes also to $0$ as $N\to\infty$. This implies \eqref{T2_claim2} and thus \eqref{T2_claim_result} which is equivalent to \eqref{T2_claim_1'}, so the proof of the theorem is now complete.

\bibliography{references}
\bibliographystyle{abbrv}

\end{document}